\documentclass[11pt]{amsart}
\usepackage{amsmath,amsthm,epsfig,amssymb}
\title{A parametrization of equilateral triangles having integer coordinates}
\author{Eugen J. Ionascu}
\curraddr{(EJI) Department of Mathematics\\ Columbus State University\\4225 University Avenue\\
Columbus, GA 31907\\
and Honorific Member of the Romanian Institute of Mathematics
``Simion Stoilow" } \email{ionascu\_eugen@colstate.edu;}
\date{\today}
\begin{document}
\def\sms{\small\scshape}
\baselineskip18pt
\newtheorem{theorem}{\hspace{\parindent}
T{\scriptsize HEOREM}}[section]
\newtheorem{proposition}[theorem]
{\hspace{\parindent }P{\scriptsize ROPOSITION}}
\newtheorem{corollary}[theorem]
{\hspace{\parindent }C{\scriptsize OROLLARY}}
\newtheorem{lemma}[theorem]
{\hspace{\parindent }L{\scriptsize EMMA}}
\newtheorem{definition}[theorem]
{\hspace{\parindent }D{\scriptsize EFINITION}}
\newtheorem{problem}[theorem]
{\hspace{\parindent }P{\scriptsize ROBLEM}}
\newtheorem{conjecture}[theorem]
{\hspace{\parindent }C{\scriptsize ONJECTURE}}
\newtheorem{example}[theorem]
{\hspace{\parindent }E{\scriptsize XAMPLE}}
\newtheorem{remark}[theorem]
{\hspace{\parindent }R{\scriptsize EMARK}}
\renewcommand{\thetheorem}{\arabic{section}.\arabic{theorem}}
\renewcommand{\theenumi}{(\roman{enumi})}
\renewcommand{\labelenumi}{\theenumi}
\newcommand{\Q}{{\mathbb Q}}
\newcommand{\Z}{{\mathbb Z}}
\newcommand{\N}{{\mathbb N}}
\newcommand{\C}{{\mathbb C}}
\newcommand{\R}{{\mathbb R}}
\newcommand{\F}{{\mathbb F}}
\newcommand{\K}{{\mathbb K}}
\newcommand{\D}{{\mathbb D}}
\def\phi{\varphi}
\def\ra{\rightarrow}
\def\sd{\bigtriangledown}
\def\ac{\mathaccent94}
\def\wi{\sim}
\def\wt{\widetilde}
\def\bb#1{{\Bbb#1}}
\def\bs{\backslash}
\def\cal{\mathcal}
\def\ca#1{{\cal#1}}
\def\Bbb#1{\bf#1}
\def\blacksquare{{\ \vrule height7pt width7pt depth0pt}}
\def\bsq{\blacksquare}
\def\proof{\hspace{\parindent}{P{\scriptsize ROOF}}}
\def\pofthe{P{\scriptsize ROOF OF}
T{\scriptsize HEOREM}\  }
\def\pofle{\hspace{\parindent}P{\scriptsize ROOF OF}
L{\scriptsize EMMA}\  }
\def\pofcor{\hspace{\parindent}P{\scriptsize ROOF OF}
C{\scriptsize ROLLARY}\  }
\def\pofpro{\hspace{\parindent}P{\scriptsize ROOF OF}
P{\scriptsize ROPOSITION}\  }
\def\n{\noindent}
\def\wh{\widehat}
\def\eproof{$\hfill\bsq$\par}
\def\ds{\displaystyle}
\def\du{\overset{\text {\bf .}}{\cup}}
\def\Du{\overset{\text {\bf .}}{\bigcup}}
\def\b{$\blacklozenge$}

\begin{abstract} We study the existence of equilateral triangles of given side
lengths and with integer coordinates in dimension three. We show
that such a triangle exists if and only if their side lengths are of
the form $\sqrt{2(m^2-mn+n^2)}$ for some integers $m,n$. We also
show a similar characterization for the sides of a regular
tetrahedron in $\Z^3$: such a tetrahedron exists if and only if the
sides are of the form $k\sqrt{2}$, for some $k\in\N$. The
classification of all the equilateral triangles in $\Z^3$ contained
in a given plane is studied and the beginning analysis is presented.
A more general parametrization is proven under a special assumption.
Some related questions are left in the end.
\end{abstract}
\maketitle
\section{INTRODUCTION}
It is known that there is no equilateral triangle whose vertices
have integer coordinates in the plane. One can easily see this by
calculating the area of such a triangle of side length $l$ using the
formula $\text{Area}=\frac{l^2\sqrt{3}}{4}$ and by using Pick's
theorem for the area of a polygon with vertices of integer
coordinates: $\text{Area}=\frac{\sharp b}{2}+\sharp i-1$ where
$\sharp b$ is the number of points of integer coordinates on the
boundary  of the polygon and $\sharp i$ is the number of such points
in the interior of the polygon. Since Pick's theorem implies that
this area is a rational number of square units, the formula
$\text{Area}=\frac{l^2\sqrt{3}}{4}$ says that this is a rational
multiple of $\sqrt{3}$ since $l^2$ must be a positive integer by the
Pythagorean theorem. This contradiction implies that no such
triangle exists.

The analog of this fact in three dimensions is not true since one
can form a regular tetrahedron by taking as vertices the points
$O(0,0,0)$, $A(1,1,0)$, $B(1,0,1)$ and $C(0,1,1)$. It turns out that
the sides of such regular tetrahedra have to be of the form
$k\sqrt{2}$, $k\in \Bbb N$. Moreover and as a curiosity, one can use
the facts derived in this note to show that there are only three
regular tetrahedrons in $\Z^3$ having the origin as one of their
vertices and of side lengths $9\sqrt{2}$, where the counting has
been done up to symmetries of the cube that leave the origin fixed.
In the figure below, that we generated with Maple, we show three
regular tetrahedrons that together with all their cube symmetries
fill out the class just described.

\begin{center}\label{figure1}
\epsfig{file=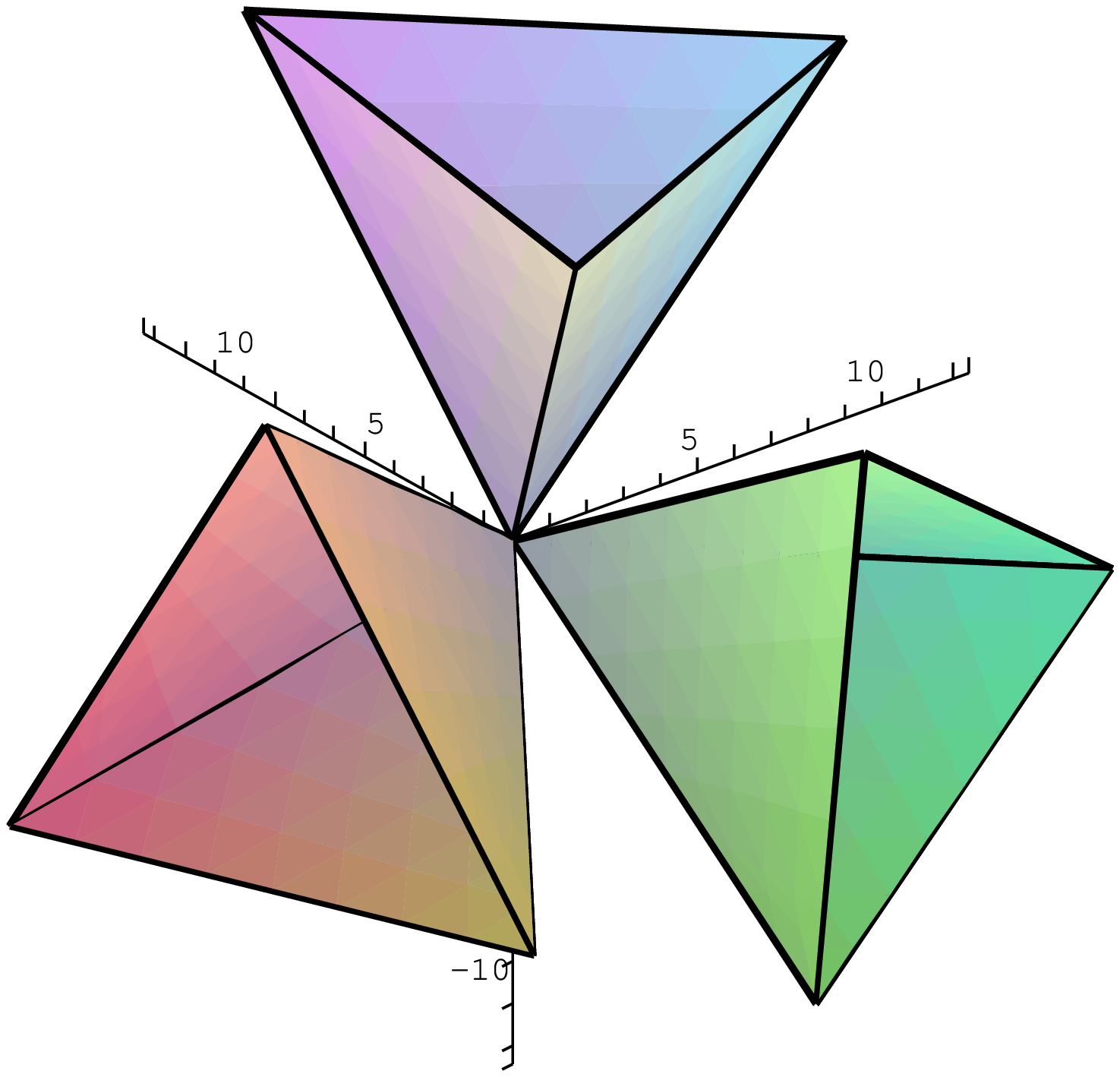,height=3in,width=3in}
\vspace{.1in}\\
{\it Figure 1: Regular tetrahedra of side lengths $9\sqrt{2}$}
\vspace{.1in}

\begin{table}\label{coordoft}
\centerline{
\begin{tabular}{|c|c|c|c|}
 \hline \hline
$O(0,0,0)$& $A_1(9,9,0)$ & $B_1(9,0,9)$& $C_1(0,9,9)$ \\
\hline
$O(0,0,0)$& $A_2(-9,9,0)$& $B_2(-4,5,-11)$& $C_2(3,12,-3)$ \\
\hline
$O(0,0,0)$& $A_3(12,3,-3)$& $B_3(7,-8,-7)$& $C_3(3,3,-12)$ \\
\hline
\end{tabular}}
\vspace{0.1in} \caption{Coordinates of the three regular tetrahedra
in Figure~\ref{figure1}}
\end{table}
\end{center}

\vspace{0.2in}

Along these lines and in retrospect to our work here,  the following
related result of Schoenberg, \cite{s}, who proved that a regular
$n$-simplex exists in $\Z^n$ in the following cases and no
others:\par

 (i) $n$ is even and $n+1$ is a square;\par

 (ii) $n\equiv 3 \ (mod\ 4)$;\par

(iii) $n\equiv  1$ $(mod\ 4)$ and $n+1$ is the sum of two
squares,\par \n seems to open a lot more questions such as: what
will be the similar corresponding parameterizations in all these
cases when regular $n$-simplexes exist? Also, one may ask when a
regular $n$-simplex exists in $\Z^m$.

Equilateral triangles with vertices of integer coordinates in the
three dimensional space are numerous as one could imagine from the
situation just described. One less obvious example is the triangle
$CDO$ with $C(31,19,76)$ and $D(44,71,11)$ having side lengths equal
to $13\sqrt{42}$. Generating all such triangles is a natural
question and we may start with one such triangle and applying the
group of affine transformations $T(\overset{\rightarrow}{x})=\alpha
O(\overset{\rightarrow}{x})+\overset{\rightarrow}{y}$ where $\alpha
\in \Q$, $O$ is an orthogonal matrix with rational coefficients and
$\overset{\rightarrow}{y}$ a vector in $\Q^3$. Such a transformation
multiplies the side lengths with the factor $\alpha$ and so for
instance the triangle $OAB$ with side lengths $\sqrt{2}$ cannot be
transformed into the $CDO$.

In fact, in the next section we show that the side lengths which
appear from such equilateral triangles are of the form
$\sqrt{2N(\zeta)}$ where $\zeta$ is an Eisenstein-Jacobi integer and
$N(\zeta)$ is its norm. The Eisenstein-Jacobi integers are defined
as $\Z[\omega]$ where $\omega$ is a primitive root of unity, i.e.
the complex numbers of the form $\zeta=m+n\omega$ with $m,n\in \Bbb
Z$ with their norm is given by $N(\zeta)=m^2+mn+n^2$.

What makes the existence of such triangles work in space that
doesn't work in two dimensions? We show in
Proposition~\ref{planeprop} that the plane containing such a
triangle must have a normal vector
$\overset{\rightarrow}{n}=a\overset{\rightarrow}{i}+b\overset{\rightarrow}{j}+c\overset{\rightarrow}{k}$
where $a,b,c$ are integers that satisfy the Diophantine equation

\begin{equation}\label{eq1}
a^2+b^2+c^2=3d^2, \ d\in \Z . \end{equation}

This equation has no non-trivial solutions, if $c=0$ for instance,
according to the Gauss' characterization for the numbers that can be
written as sums of two perfect squares.

Our study of the existence of such triangles started with an
American Mathematics Competition problem in the beginning of 2005.
The problem was:
\begin{problem}\label{originalpr} Determine the number of equilateral triangles whose
vertices have coordinates in the set $\{0,1,2\}$.
\end{problem}
It turns out that the answer to this question is 80.  Let us
introduce the notation $\cal E\cal T(n)$ for the number of
equilateral triangles whose vertices have coordinates in the set
$\{0,1,2, ...,n\}$ for $n\in \N$. Some of the values of $\cal E\cal
T(n)$  are tabulated next:\par

\vspace{0.1in}
\begin{table}\label{nofeqtrinagles0n}
 \centerline{
\begin{tabular}{|c|c|c|c|c|c|c|c|c|c|c|}
  \hline
   $n$ & $1$ & $2$ & $3$& $4$& $5$& $6$& $7$& $8$ & $9$ & $10$  \\  \hline
$\cal E\cal T(n)$& $8$ & $80$ &  $368$& $1264$ & $3448$& $7792$& $16176$& $30696$ & $54216$& $90104$ \\
\hline
\end{tabular}}
\vspace{0.1in} \caption{Sequence A 102698}
\end{table}

\vspace{0.1in} This sequence has been entered in the on-line
Encyclopedia of Integers Sequences by Joshua Zucker on February 4th,
2005. The first $34$ terms in this sequence were calculated by Hugo
Pfoertner using a program in Fortran. Our hope is that the results
obtained here may be used in designing a program that could
calculate $\cal E\cal T(n)$ for significantly more values of $n$.

\section{Planes containing equilateral triangles in $\Z^3$}
Let us denote the side lengths of an equilateral triangle $\triangle
OPQ$ by $l$. We are going to discard translations, so we may assume
that one of the vertices of such a triangle is $O(0,0,0)$. If the
other two points, $P$ and $Q$, have coordinates $(x,y,z)$ and
$(u,v,w)$ respectively, then as we have seen before the area of
$\triangle OPQ$ is given by

\begin{equation}\label{area}
\text{Area}=\frac{l^2\sqrt{3}}{4}=\frac{1}{2}|\overset{\rightarrow}{OP}
\times \overset{\rightarrow}{OQ}|=\left|
                                    \begin{array}{ccc}
                                      \overset{\rightarrow}{i} & \overset{\rightarrow}{j} & \overset{\rightarrow}{k} \\
                                      x & y & z \\
                                      u & v & w \\
                                    \end{array}
                                  \right|
.\end{equation}

This implies the following simple proposition but essential in our
discussion:

\begin{proposition}\label{planeprop}
Assume the triangle $\triangle OPQ$ is equilateral and its vertices
have integer coordinates with $O$ the origin and $l=|OP|$. Then the
points $P$ and $Q$ are contained in a plane of equation
$a\alpha+b\beta+c\gamma =0$, where $a,b,c$, and $d$ are integers
which satisfy (\ref{eq1}) and $l^2=2d$.
\end{proposition}
\proof. \ Assume the coordinates of $P$ and $Q$ are denoted as
before. Let us observe that
$(x-u)^2+(y-v)^2+(z-w)^2=2l^2-2(xu+yv+zw)$ which implies
$xu+yv+zw=\frac{l^2}{2}=d\in \Z$. Then using the fact that
$\overset{\rightarrow}{OP}$ and $\overset{\rightarrow}{OQ}$ are
contained in the plane orthogonal on the vector
$\overset{\rightarrow}{OP} \times
\overset{\rightarrow}{OQ}=a\overset{\rightarrow}{i}+b\overset{\rightarrow}{j}+c\overset{\rightarrow}{k}$
with $a=yw-vz$, $b=zu-xw$ and $c=xv-yu$ the statement follows from
(\ref{area}).\eproof

The equation (\ref{eq1}) has infinitely many integer solutions
besides the obvious ones $a=\pm d$, $b=\pm d$, $c=\pm d$. For
instance we can take $a=-19$, $b=11$, $c=5$ and $d=13$ and the
triangle $OCD$ given in the Introduction has $C$ and $D$ in the
plane $\{(\alpha,\beta,\gamma)\in \R^3|\ -19\alpha+11\beta+5\gamma
=0\}$.
\begin{definition}\label{planesig}
Let us introduce the set $\cal P$ of all these planes, i.e. $\{
(\alpha,\beta,\gamma)\in \R^3|\ a\alpha+b\beta+c\gamma =e\} $, such
that $a^2+b^2+c^2=3d^2$ for some $a,b,c,d,e\in \Z$ and
$gcd(a,b,c)=1$.
\end{definition}

So, if one starts with a plane $\pi$ in $\cal P$, picks two points
of integer coordinates that belong to $\pi$, which amounts to
solving a simple linear Diophantine equation, the natural question
is weather or not there exists a third point of integer coordinates,
contained in $\pi$, that completes the picture to an equilateral
triangle (see Figure~\ref{figure2}). In order for the third point to
exist one needs to take the first two points in a certain way. But
if one requires only that the new point have rational coordinates it
turns out that this is always possible and the next theorem gives a
way to find the coordinates of the third point in terms of the given
data.

\begin{theorem}\label{construction}
Assume that $P(u,v,w)$ ($u,v,w\in\Q$) is an arbitrary point of a
plane $\pi\in \cal P$ of normal vector $(a,b,c)$  and passing
through the origin $O$. Then the coordinates of a point $Q(x,y,z)$
situated in $\pi$ and such that the triangle $\triangle OPQ$ is
equilateral are all rational numbers given by:
\begin{equation}\label{eqfortp}
\begin{cases}
\ds x=\frac{u}{2}\pm \frac{cv-bw}{2d}\\
\ds y=\frac{v}{2}\pm \frac{aw -cu}{2d}  \\
\ds z =\frac{w}{2}\pm \frac{bu-av}{2d}.
\end{cases}
\end{equation}
\end{theorem}

\begin{center}\label{figure2}
\epsfig{file=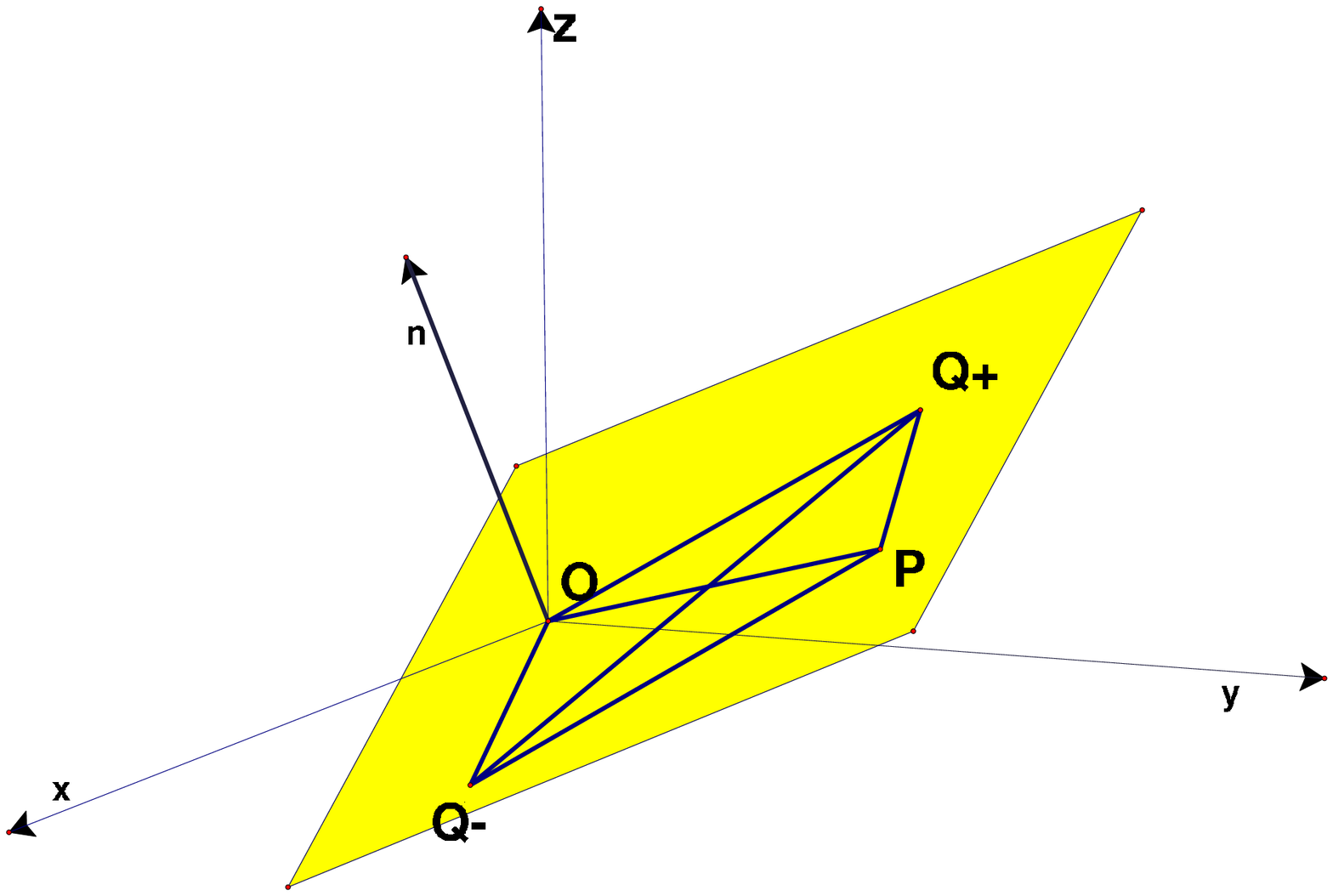,height=2in,width=3in}
\vspace{.1in}\\
{\it Figure 2: Plane of normal (a,b,c)} \vspace{.1in}
\end{center}

\vspace{0.2in} \proof. \ From the geometric interpretation of the
problem of finding a point $Q\in \pi$ such that $\triangle OQP$
becomes equilateral one can see that there are only two points that
satisfy this condition and just having real coordinates. We want to
show that one point is given by taking the plus sign in all
equalities in (\ref{eqfortp}) and the other point corresponds to the
minus sign in (\ref{eqfortp}). We are going to set
$\overset{\rightarrow}{n}=\ds \frac{1}{d\sqrt{3}}(a,b,c)$ which is
one of the two unit vectors normal to the plane $\pi$ and let
$\overset{\rightarrow}{r}=\overset{\longrightarrow}{OP}=(u,v,w)$.
Then the cross product $\overset{\rightarrow}{r}\times
\overset{\rightarrow}{n}$ is given by

$$\overset{\rightarrow}{r}\times \overset{\rightarrow}{n}=\ds \frac{1}{d\sqrt{3}}(cv-bw,aw -cu,bu-av).$$
So we observe that the solution $(x,y,z)$ is, in fact, if written in
vector notation,
$\overset{\longrightarrow}{OQ}_{\pm}=\frac{1}{2}\overset{\rightarrow}{r}\pm
\frac{\sqrt{3}}{2}\overset{\rightarrow}{r}\times
\overset{\rightarrow}{n}$. It is easy now to check that
$\overset{\rightarrow}{r}$ and $\overset{\longrightarrow}{OQ}_{\pm}$
have the same norm and make a $60^{\circ}$ angle in between.\eproof

\section{Solutions of the Diophantine equation $a^2+b^2+c^2=3d^2$}

We have enough evidence to believe that for each $(a,b,c)$
satisfying (\ref{eq1}) for some $d\in \Z$, there are infinitely many
equilateral triangles in $\Z^3$ that belong to a plane of normal
vector $(a,b,c)$ and the purpose of this paper is to determine a way
to generate all these triangles. But how many planes do we have in
$\cal P$?  Let us observe that if $d$ is even, then not all of $a$,
$b$, $c$ can be odd integers, so at least one of them must be even.
Then the sum of the other two is a number divisible by $4$ which is
possible only if they are also even. Therefore we may reduce all
numbers by a factor of two in this case. So, if we assume without
loss of generality that $gcd(a,b,c)=1$, then such solutions of
(\ref{eq1}) must have $d$ and odd integer and then this forces that
$a$, $b$ and $c$ must be all odd integers too. If in addition,
disregarding the signs, we have $a,b,c\in\N$ and $a\le b\le c$, then
such a solution will be referred to as a {\it primitive} solution of
(\ref{eq1}). One can find lots of solutions of (\ref{eq1}) in the
following way.

\begin{proposition}\label{paramofplane} (i) The following formulae give a three integer
parameter solution of (\ref{eq1}):
\begin{equation}\label{peqabc}
\begin{cases}
a=-x_1^2+x_2^2+x_3^2-2x_1x_2-2x_1x_3\\
b=x_1^2-x_2^2+x_3^2-2x_xx_1-2x_2x_3\\
c=x_1^2+x_2^2-x_3^2-2x_3x_1-2x_3x_2\\
d=x_1^2+x_2^2+x_3^2
\end{cases},\ \  x_1,x_2,x_3\in \Z.
\end{equation}
\par

(ii) Every primitive solution of (\ref{eq1}) is of the form
(\ref{peqabc}) with $x_1,x_2,x_3 \in \Q \sqrt{k}$ with some $k\in
\N$.

\end{proposition}
\n \proof.\  (i) Let us observe that $a=d-2x_1s$, $b=d-2x_2s$,
$x_3=d-2x_3s$ where $s=x_1+x_2+x_3$. Then a simple calculation shows
that $a^2+b^2+c^2=3d^2$. In fact, for every $x_1,x_2,x_3\in \Q$ for
which $a,b,c,d$ given by (\ref{peqabc}) are integers we get a
solution of (\ref{eq1}).

If one has a solution $(a,b,c,d)$ of (\ref{eq1}) with $d\ne 0$, then
we can introduce $x_1=(d-a)t$, $x_2=(d-b)t$, $x_3=(d-c)t$ where
$t\in \R$ such that $x_1^2+x_2^2+x_3^2=d$. This gives
$$t^2=\frac{d}{(d-a)^2+(d-b)^2+(d-c)^2}=\frac{d}{6d^2-2d(a+b+c)}=
\frac{1}{2(3d-a-b-c)}=\frac{1}{2}\frac{t}{x_1+x_2+x_3},$$

\n and so  $t=\frac{1}{2s}$ with
$s=x_1+x_2+x_3=\sqrt{(3d-a-b-c)/2}$. Every primitive solution of
(\ref{eq1}) must have $a$, $b$, $c$, and $d$ all odd numbers as we
have observed. This make $k=(3d-a-b-c)/2$ an integer. Since $3d\ge
a+b+c$ is equivalent to
$3(a^2+b^2+c^2)-(a+b+c)^2=(a-b)^2+(b-c)^2+(a-c)^2\ge 0$, it follows
that $k\in \N$. Therefore in general every solution is of the form
in (\ref{peqabc}) but with $x_1,x_2,x_3 \in \Q \sqrt{k}$.\eproof

\vspace{0.1in} \n  {\bf Example:} Suppose $a=5$, $b=11$, $c=19$ and
$d=13$ then $x_1=2\sqrt{2}$, $x_2=\frac{1}{2}\sqrt{2}$ and
$x_2=-\frac{3}{2}\sqrt{2}$.

\vspace{0.1in} \n Although formulae (\ref{peqabc}) may provide
infinitely many solutions of (\ref{eq1}), the following proposition
is giving more information about its integer solutions.

\begin{proposition}\label{sofabcfedodd} The equation (\ref{eq1}) has non-trivial
solutions for every odd integer $d\ge 3$.
\end{proposition}

\n \proof. \   If $d=2p+1$ for some $p\in \Z$, $p\ge 1$, then
$3d^2=3[4p(p+1)+1]= 8l+3$ which shows that $3q^2\equiv 3$  (mod 8).
From Gauss's Theorem about the number of representations of a number
as a sum of three squares (Theorem 2, pp. 51 in \cite{gr}), we see
that the number of representations of $3d^2$ as a sum of three
squares is at least 24. Here, the change of signs is counted so each
solution actually generates eight solutions by just changing the
signs. Also, the six permutations of $a$, $b$ and $c$ get into the
counting process. So, there must be at least one solution which is
nontrivial since the solution $3d^2=d^2+d^2+d^2$ generate only eight
solutions by using all possible change of signs and all
permutations.~\eproof

Another way of generating an infinite family of solutions of
(\ref{eq1}) is to reduce it to two separate equations, say:
$146=3d^2-c^2$ and $ a^2+b^2=146$. The second equation admits as
solution, for instance, $a=11$ and $b=5$. The first equation has a
particular solution $d=7$ and $c=1$. Then one can use the recurrence
formulae to obtain infinitely many solutions of $146=3d^2-c^2$:
$$d_{n+1}=2d_n+c_n, c_{n+1}=3d_n+2c_n\  for\  n\in \Bbb N$$ and
$d_1=7$, $c_1=1$. A simple calculation shows that
$3d_{n+1}^2-c_{n+1}^2=3d_n^2-c_n^2$ so, by induction, $(d_n, c_n)$
is a solution of the equation $3d^2-c^2=146$ for all $n\in \Bbb N$.
It is easy to see that $q_n$ and $c_n$ are increasing sequences so
this procedure generates infinitely many solutions of (\ref{eq1}).

Some of the primitive solutions of (\ref{eq1}) are included in the
table below.

\vspace{0.1in}
\begin{table}\label{psoabc}
\centerline{
\begin{tabular}{||c|c||c|c||}
 \hline \hline
 d& (a,b,c) &  d & (a,b,c )\\
 \hline \hline
$1$& $\{(1,1,1)\}$ & $9$& $\{(1,11,11),(5,7,13)\}$\\
\hline $3$& $\{(1,1,5)\}$ & $11$& $\{(1,1,19),(5,7,17),(5,13,13)\}$\\
\hline $5$& $\{(1,5,7)\}$& $13$& $\{(5,11,19),(7,13,17)\}$\\
\hline $7$& $\{(1,5,11)\}$ & $15$& $\{(1,7,25),(5,11,25),(5,17,19)\}$\\
\hline \hline
\end{tabular}}\vspace{0.1in}
\caption{Primitive solutions of (\ref{eq1})} \vspace{0.1in}
\end{table}
As a curiosity the number of representations of $d=1003$ (which will
give the number of ``primitive" representations of
$a^2+b^2+c^2=3\times 2006^2$) is 182.

\section{Generating the beginning of the parametrization family}

The simplest solution of (\ref{eq1}) is $a=b=c=d=1$. We are going to
introduce some more notation here before we give the parametrization
for this case.

\begin{definition}\label{etinplabc} For every $(a,b,c)$, a primitive
solution of (\ref{eq1}), denote by ${\cal T} _{a,b,c}$ the set of
all equilateral triangles with integer coordinates having the origin
as one of the vertices and the other two lie in the plane
$\{(\alpha,\beta,\gamma)\in \R^3 |a\alpha +b\beta+c\gamma=0\}$.
\end{definition}

It is clear now that in light of  Proposition~\ref{planeprop}, every
equilateral triangle having integer coordinates after a translation,
interchange of coordinates, or maybe a  change of signs of some of
the coordinates, belongs to one of the classes $\cal T_{a,b,c}$. Let
us introduce also the notation $\cal T$ for all the equilateral
triangles in $\Z^3$. If we have different values for $a$, $b$ and
$c$, how many different planes can one obtain by permuting  $a$, $b$
and $c$ in between and changing their signs? That will be $6$
permutations and essentially $4$ change of signs (note that
$ax+by+cz=0$ is the same plane as $(-a)x+(-b)y+(-c)z=0$) which gives
a total of $24$ such transformations. We are going to denote the
group of symmetries of the space determined by these transformations
and leave the origin fixed, by $\cal S_{cube}$ (it is actually the
group of symmetries of the cube). Hence we have

\vspace{0.1in}
\begin{equation}\label{partition}
{\cal T} =\underset{\begin{array}{c}  s\in {\cal S}_{cube},s(O)=O,\\
a^2+b^2+c^2=3d^2\\ 0<a\le b\le c, gcd(a,b,c)=1\\a,b,c,d\in \Z , v\in
\Z^3 \end{array} }{ \bigcup }  s \left ({\cal T}_{a,b,c}\right )+v.
\end{equation}

\vspace{0.1in}

\begin{theorem}\label{firstpar} Every triangle $OAB \in \cal T_{1,1,1}$ is of the
form $\{A,B,O\}=\{(m, -n,n-m), (m-n,-m,n), (0,0,0)\}$ for some
$m,n\in \Bbb Z$. The side lengths of $\triangle OAB_{m,n}$ are given
by $$l=\sqrt{2(m^2-mn+n^2)}.$$
\end{theorem}
\n \proof. \  Let us assume $A$ has coordinates $(u,v,w)$ with
$u+v+w=0$ and $B(x,y,z)$ with $x+y+z=0$. From (\ref{eqfortp}) we get
that $x=\frac{u}{2}+\frac{v-w}{2}$, $y=\frac{v}{2}+\frac{w-u}{2}$,
$z=\frac{v}{2}+\frac{u-v}{2}$ if we choose the plus signs. This
choice is without loss of generality since we can interchange the
roles of $A$ and $B$ if necessary. This implies $x=-w$, $y=-u$ and
$z=-v$. So, if we denote $u=m$, $v=-n$ then $w=n-m$ and so $x=m+n$,
$y=-m$, $z=-n$.\eproof

We are introducing the notation $N(m,n)=2(m^2-mn+n^2)$ for
$m,n\in\Z$. For the next cases, $d\in \{3,5,7\}$, as we have
recorded in the Table \ref{psoabc}, $3d^2$ has also a unique
primitive representation. One can use basically the same technique
as in the proof of Theorem~\ref{firstpar} to derive the
corresponding parameterizations for the vertices in $\cal T_{a,b,c}$
($d\in \{3,5,7\}$) and the corresponding side lengths but for each
individual set of formulae, that are given below, we had something
specific to work out in order to get rid of denominators that
naturally arise from applying (\ref{eqfortp}):

\[
\begin{array}{c}
d=3,  \ l=3\sqrt{N(m,n)} \\ \cal T_{1,1,5}=\{[O,
(4m-3n,m+3n,-m),(3m+n,-3m+4n,-n)]:m,n\in \Z, l\neq 0 \},
\\ \\
d=5, \ l=5\sqrt{N(m,n)} \\ \cal T_{1,5,7}=\{[O,
(7m-4n,5n,-m-3n),(3m-7n,5m,-4m+n)]:m,n\in \Z , l\neq 0 \},
\\ \\
d=7, \ l=7\sqrt{N(m,n)} \\ \cal T_{1,5,11}=\{[O,
(8m-9n,5m+4n,-3m-n),(-m-8n,9m-5n,-4m+3n)]:m,n\in \Z, l\neq 0\}.
\end{array}
\]
\vspace{0.1in}

{\bf Remark:} Every triangle in one particular family, $s(\cal
T_{a,b,c})+v$, is different of all the other triangles in other
families since they live in different planes. So if we take in
(\ref{partition}) only the symmetries $s\in {\cal S}_{cube}$,
$s(O)=O$, that give different normal vectors (two of the numbers
$a$, $b$, $c$ may be equal) then (\ref{partition}) is a partition of
$\cal T$.

The case $d=9$ is the first in which there are two essentially
different primitive representations of $3d^2$: $3 \times
9^2=1^2+11^2+11^2+1=5^2+7^2+13^2$. The corresponding
parameterizations are included here for completion:

\[
\begin{array}{c}
d=9,  \ l=9\sqrt{N(m,n)}\\  \\ \cal T_{1,11,11}=\{[O,
(11m-11n,4m+5n,-5m-4n),(-11n,9m-4n,-9m+5n)]:m,n\in \Z, l\neq 0 \},
\\ \\ \cal T_{5,7,13}=\{[O,
(7m+5n,8m-11n,-7m+4n),(12m-7n,-3m-8n,-3m+7n)]:m,n\in \Z , l\neq 0
\}.
\end{array}
\]
\vspace{0.1in}

To give and idea of how we obtained these parametrizations we will
include the proof of the case $d=9$, $a=5$, $b=7$, $c=13$, that gave
us  the last of the above formulae.

\n {\bf \proof.} \  Assume that one of the points, $P$,  has
coordinates $(u,v,w)$. If one solves the Diophantine equation
$5u+7v+13w=0$ finds that a general solution may be written as
\[
\begin{cases}
w=5u+7t,\\
v=-10u-13t \ and\  t,w\in \Z.
\end{cases}
\]

 Using Theorem~\ref{construction} we see
that the coordinates of a point $Q$, say $(x,y,z)$, such that
$\triangle OPQ\in \cal T_{5,7,13}$ must be given by (\ref{eqfortp}).
Switching $P$ with $Q$, if necessary, we may assume that the signs
in (\ref{eqfortp}) could be taken all plus signs. This gives
\[
\begin{cases}
\ds x=-\frac{26u}{3}-\frac{109t}{9},\\ \\
\ds y=-\frac{13u}{3}-\frac{41t}{9},\\ \\
\ds z=\frac{17u}{3}+\frac{64t}{9}.
\end{cases}
\]

\n Since $x=-9u-12t-\frac{t-3u}{9}$  must be an integer we need to
have $t=3u+9g$ for some $g\in \Z$. Substituting we find  that all
other coordinates are integers: $x=-45u-109g$, $y=-18u-41g$,
$z=27u+64g$. Calculating $l^2=u^2+v^2+w^2$ we get $l^2=
3078u^2+14742ug+17658g^2=2\times 9^2(19u^2+91ug+109g^2)$. Or
$l^2=2\times 9^2(19u^2+91ug+109g^2)=2\times
9^2[(2u+5g)^2+(2u+5g)(3u+7g)+(3u+7g)^2]$ which suggests that we
could change the variables $2u+5g=-m$ and $3u+7g=n$ to obtain the
statement from the theorem. If we solve this system it turns out
that the solution preserves integers values since $u=7m+5n$ and
$g=-3m-2n$. \eproof

A natural question that we may ask at this point is whether or not
every $\cal T_{a,b,c}$ admits such a parametrization. In the next
section we prove that this is indeed the case under the assumption
that $\min \{gcd(a,d),gcd(b,d),gcd(c,d)\}=1$. We checked all the odd
integers $d$ between  $1$ and $401$ and this condition is satisfied
for all such $d$'s.

\section{Characterization of side lengths}

We will begin with two preliminary results.  The first we just need
to recall it since it is a known fact that can be found in number
theory books mostly as an exercise or as an implicit corollary of
more general theorems about quadratic forms or Euler's $6k+1$
theorem (see \cite{r}, pp. 568 and \cite{g}, pp. 56).
\begin{proposition}\label{ejtypet} An integer $t$ can be written as $m^2-mn+n^2$ for some
$m,n\in \Z$ if and only if in the prime factorization of $t$, $2$
and the primes of the form $6k-1$ appear to an even exponent.
\end{proposition}


The next lemma is probably also known in algebraic number theory but
we do not have straight reference for it so we are going to include
a proof of it.

\begin{lemma}\label{thebeauty}  An integer $t$ which can be written as  $t=3x^2-y^2$ with
$x,y \in \Bbb Z$ is the sum of two squares if and only if $t$ is of
the form $t=2(m^2-mn+n^2)$ for some integers $m$ and $n$.
\end{lemma}

\proof . \ For necessity, by Proposition~\ref{ejtypet}, we have to
show that $t$ is even and $t/2$ does not contain in its prime factor
decomposition any of the primes $2$ or those of the form $6k-1$
except to an even power. First, let us show that $t$ must be even
and the exponent of $2$ in its prime factorization is odd. Since
$t=3x^2-y^2=a^2+b^2$ implies $3x^2=a^2+b^2+y^2$ we have observed
that either all $x$, $y$, $a$, and $b$ are even or all odd.

If $x$, $y$, $a$ and $b$ are all even we can factor out a $2$ from
all these numbers and reduce the problem to $t/4$ instead of $t$.
Applying this arguments several times one can see that
$t=2^{2l+1}t'$ with $t'$ odd and $l\in \Bbb Z$. Without loss of
generality we may assume that $l=0$. In this case $t$ contains only
one power of $2$ in its prime decomposition and so $x$, $y$, $a$ and
$b$ must be all odd.

Let us then suppose that $t/2$ is divisible by a prime $p=6k-1$ for
some  $k\in \Bbb N$. We need to show that the exponent of $p$ in the
prime factorization of $t$ is even. Since $p$ divides $t=3x^2-y^2$
we get that $3x^2\equiv y^2$ (mod p). If $p$ divides $x$, then $p$
divides $y$ and so $p^2$ divides $t$ which reduces the problem to
$t/p^2$. Applying this argument several times we arrive to a point
when $p$ does not divide $x/p^i$. So, discarding an even number of
$p$'s from $t$, we may assume that $i=0$. This implies that $x$ has
an inverse modulo $p$ and then $z^2\equiv 3$ (mod $p$), where
$z=x^{-1}y$. Using the Legendre symbol this says that
$(\frac{3}{p})=1$. By the Law of Quadratic Reciprocity (Theorem
11.7, in \cite{r}) we see that
$(\frac{p}{3})=(-1)^{\frac{(p-1)}{2}\frac{3-1}{2}}=(-1)^{3k-1}$. But
the equation in $w$, $w^2\equiv p$ (mod $3$), is equivalent to
$w^2\equiv -1$ (mod $3$) which obviously has no solution in $w$.
This implies $(\frac{p}{3})=-1$ and so $k$ has to be even. Therefore
$p=12k'-1=4j+3$ for some $j\in \Bbb Z$. But by hypothesis, $t$ is a
sum of two squares and so, from Euler's characterization of those
numbers, $p$ must have an even exponent in the prime decomposition
of $t$.

For sufficiency, let us assume that $t=3x^2-y^2=2(m^2-mn+n^2)$ for
some $x,y,m,n\in \Bbb Z$. Using Euler's characterization we have to
show that if $p=4k+3$ is a prime dividing $t$ then the exponent in
its prime decomposition is even. If $p=3$, then $3$ divides $y$
which implies $x^2-3{y'}^2=2({m'}^2-m'n'+{n'}^2)$. This is true
because of Proposition~\ref{ejtypet} which one has to use in both
directions.

If ${m'}^2-m'n'+{n'}^2$ is not divisible by $3$ then
${m'}^2+m'n'+{n'}^2\equiv 1$ (mod 3) since all the prime factors of
the form $6k-1$ and $2$ appear to even exponents. This implies
$x^2-3{y'}^2\equiv 2({m'}^2-m'n'+{n'}^2)$ (mod $3$) or $x^2\equiv 2$
(mod $3$) which is a contradiction. So, ${m'}^2-m'n'+{n'}^2$ must
contain another factor of $3$ and so the problem could be then
reduced to $t/9$ instead of $t$. Hence $3$ must have an even
exponent in the prime decomposition of $t$.

Let us assume that $k=3j-1$ with $j\in \Bbb Z$. Then
$p=12j-1=6(2j)-1$ and so these primes must appear to an even power
in the decomposition of $m^2-mn+n^2$. The case $k=3j+1$ ($p=12k+7$)
is not possible because that will contradict the Law of Quadratic
Reciprocity: $(\frac{3}{p})(\frac{p}{3})=1\not
=(-1)^{\frac{(p-1)}{2}\frac{3-1}{2}}$. \eproof

\begin{theorem}\label{main}
An equilateral triangle of side lengths $l$ and having integer
coordinates in $\Bbb R^3$ exists, if and only if
$l=\sqrt{2(m^2-mn+n^2)}$ for some integers $m$ and $n$ (not both
zero) .
\end{theorem}
\n {\bf \proof.} \ The sufficiency part of the theorem is given by
the triangles in $\cal T_{1,1,1}$ (Theorem~\ref{firstpar}). For
necessity let us start with an arbitrary equilateral triangle having
integer coordinates and non zero side lengths $l$. Without loss of
generality we may assume that one of its vertices is the origin.
Denote the triangle as before $\triangle OPQ$, with $P(u,v,w)$ and
$Q(x,y,z)$. As we have shown in Proposition~\ref{planeprop} we know
that $au+bv+wc=0$ and $ax+by+cz=0$ for some $a,b,c$ satisfying
$a^2+b^2+c^2=3d^2$ and $gcd(a,b,c)=1$. We noticed too that all
$a,b,c$ have to be odd integers and so, in particular, they are all
non-zero numbers.

Then $$l^2=u^2+v^2+w^2=\left(\ds
\frac{bv+cw}{a}\right)^2+v^2+w^2=\ds
\frac{(a^2+b^2)v^2+2bcvw+(a^2+c^2)w^2}{a^2}.$$

Completing the square we have

\[
\begin{array}{c}
a^2l^2=(a^2+b^2)v^2+2bcvw+(a^2+c^2)w^2= (3d^2-c^2)\left (v+\ds
\frac{bcw}{3 d^2-c^2}\right )^2+(3d^2-b^2)w^2-\ds
\frac{b^2c^2w^2}{3d^2-c^2}=\\ (3d^2-c^2)\left (v+\ds
\frac{bcw}{3d^2-c^2}\right )^2+3\frac{d^2a^2w^2}{3d^2-c^2},
\end{array}\]

or

$$a^2(3d^2-c^2)l^2=[(3d^2-c^2)v+bcw]^2+3d^2a^2w^2.$$

This calculation shows that $a^2(3d^2-c^2)l^2=m'^2-m'n'+n'^2$ where
$m'=(3d^2-c^2)v+bcw+daw$ and $n'=2daw$. By Lemma~\ref{thebeauty} we
can write $3d^2-c^2=a^2+b^2=2(m''^2-m''n''+n''^2)$ for some
$m'',n''\in\Z$. Hence
$$l^2=2\frac{1}{(2a)^2}\frac{m'^2-m'n'+n'^2}{m''^2-m''n''+n''^2}.$$ Because
$l^2\in \Z$ and Proposition~\ref{ejtypet} we see that
$l^2=2(m^2-mn+n^2)$. \eproof

We include here a similar result and the last of this section which
is only based on Proposition~\ref{planeprop}.
\begin{proposition}\label{tetrahedra}
A regular tetrahedra of side lengths $l$ and having integer
coordinates in $\Bbb R^3$ exists, if and only if $l=m\sqrt{2}$ for
some $m\in \N$.
\end{proposition}

\n {\bf \proof.} For sufficiency, we can take the tetrahedra $OPQR$
with $P(m,0,m)$, $Q(m,m,0)$ and $R(0,m,m)$.

For necessity, without loss of generality we assume the tetrahedra
$OPQR$ is regular and has all its coordinates integers. As before,
we assume $P(u,v,w)$ and $Q(x,y,z)$.

\begin{center}\label{figure3}
\epsfig{file=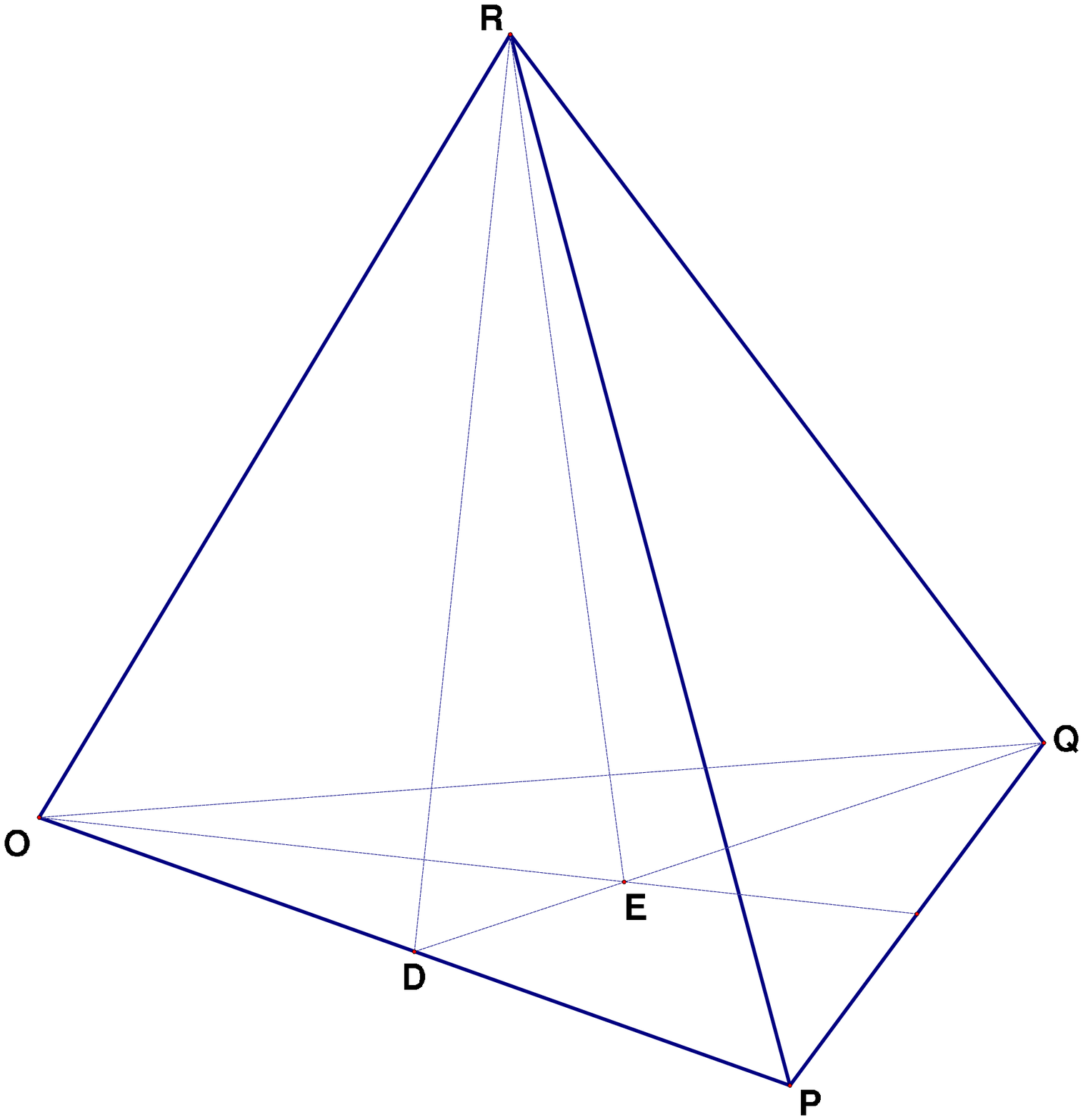,height=3in,width=3in}
\vspace{.1in}\\
{\it Figure 3: Regular tetrahedra} \vspace{.1in}
\end{center}

Let $E$ be the center of the face $\triangle OPQ$. Then from
Proposition~\ref{planeprop} we know that
$\frac{\overset{\rightarrow}{ER}}{|\overset{\rightarrow}{ER}| }=
\frac{(a,b,c)}{\sqrt{a^2+b^2+c^2}}=\frac{1}{\sqrt{3}d}(a,b,c)$ for
some $a,b,c,d\in\Z$, $l^2=2d$. The coordinates of $E$ are
$(\frac{u+x}{3},\frac{y+v}{3},\frac{z+w}{3})$.

From the Pythagorean theorem one can find easily that $RE=l
\sqrt{\frac{2}{3}}$. Since
$\overset{\rightarrow}{OR}=\overset{\rightarrow}{OE}+\overset{\rightarrow}{ER}$,
the coordinates of $R$ must be given by

$$\left(\frac{u+x}{3}\pm
l\sqrt{\frac{2}{3}}\frac{1}{\sqrt{3}d}a,\frac{y+v}{3}\pm
l\sqrt{\frac{2}{3}} \frac{1}{\sqrt{3}d}b,\frac{z+w}{3}\pm
l\sqrt{\frac{2}{3}}\frac{1}{\sqrt{3}d}c\right)$$ or
$$\left(\frac{u+x}{3}\pm \frac{2\sqrt{2}}{3l}a,\frac{y+v}{3}\pm \frac{2\sqrt{2}}{3l}b,
\frac{z+w}{3}\pm \frac{2\sqrt{2}}{3l}c\right).$$ Since these
coordinates are assumed to be integers we see that
$l=m\sqrt{2}$.\eproof




A natural question that one may ask at this point is weather or not
every family $\cal T_{a,b,c}$ contains triangles which are faces of
regular tetrahedra with integer coordinates. We believe that every
triangle in $\cal T_{a,b,c}$ is the face of such a tetrahedra as
long as its sides, in light of Theorem~\ref{main}, are of the form
$l=\sqrt{2N(m,n)}$ with $N(m,n)$ a perfect square. We leave this
conjecture for further study. To find values of $m,n\in\ Z$ such
that $N(m,n)=m^2-mn+n^2$ is a perfect square, of course one can
accomplish this in the trivial way, by taking $m=0$ or $n=0$ but
there are also infinitely many non-trivial solutions as one can see
from Proposition~\ref{ejtypet}.

\section{A more general parametrization}

Our construction depends on a particular solution, $(r,s)\in \Z^2$,
of the equation:
\begin{equation}\label{rs}
2(a^2+b^2)=s^2+3r^2.
\end{equation}

 As before let us assume that $a$, $b$, $c$ and $q$ are integers satisfying $a^2+b^2+c^2=3d^2$ with $d$ an
odd positive integer and $gcd(a,b,c)=1$. By Lemma~\ref{thebeauty} we
see that $3d^2-c^2=a^2+b^2=2(f^2-fg+g^2)$ for some $f,g\in \Z$ and
so $2(a^2+b^2)=(2f-g)^2+3g^2$ which says that equation (\ref{rs})
has always an integer solution.

\begin{theorem}\label{generalpar} Let $a$, $b$, $c$, $d$ be odd positive
integers such that $a^2+b^2+c^2=3d^2$, $a\le b\le c$ and
$gcd(d,c)=1$. Then $\cal T_{a,b,c}=\{ \triangle OPQ|\ m,n\in\Z\}$
where the points $P(u,v,w)$ and $Q(x,y,z)$  are given by
\begin{equation}\label{paramone}
\begin{cases}
u=m_um-n_un,\\
v=m_vm-n_vn,\\
w=m_wm-n_wn, \\
\end{cases}
\ \ \ \text{and} \ \ \
\begin{cases}
x=m_xm-n_xn,\\
y=m_ym-n_yn,\\
z=m_zm-n_zn,
\end{cases}\ \
\end{equation}
with
\begin{equation}\label{paramtwo}
\begin{array}{l}
\begin{cases}
m_x=-\frac{1}{2}[db(3r+s)+ac(r-s)]/q,\ \ \ n_x=-(rac+dbs)/q\\
m_y=\frac{1}{2}[da(3r+s)-bc(r-s)]/q,\ \ \ \ \  n_y=(das-bcr)/q\\
m_z=(r-s)/2,\ \ \ \ \ \ \ \ \ \ \ \ \ \ \ \  \ \ \ \ \  \ \ \ \ \ \ n_z=r\\
\end{cases}\\
\ \ \ \text{and}\ \ \\
\begin{cases}
m_u=-(rac+dbs)/q,\ \ \ n_u=-\frac{1}{2}[db(s-3r)+ac(r+s)]/q\\
m_v=(das-rbc)/q,\ \ \ \ \   n_v=\frac{1}{2}[da(s-3r)-bc(r+s)]/q\\
m_w=r,\ \ \ \ \ \ \ \ \ \ \ \ \ \ \ \ \ \ \  n_w=(r+s)/2\\
\end{cases}
\end{array}
\end{equation}
where $q=a^2+b^2$ and $(r,s)$ is a suitable solution of (\ref{rs}).
\end{theorem}

In order to show Proposition~\ref{generalpar} we need the following
lemma.

\begin{lemma}\label{existence}
Suppose $A$ and $B$ are integers in such a way $A^2+3B^2$ is
divisible by $q$ where $2q$ can be written as $s'^2+3r'^2$ with
$r',s'\in \Z$. Then there exist a writing $2q=s^2+3r^2$, $r,s\in
\Z$, such that
$$Ar+Bs \equiv 0\ \ (mod\ 2q),$$
and
$$As-3Br \equiv 0\ \ (mod \ 2q).$$
\end{lemma}

\proof. Let us observe that every number of the form $A^2+3B^2$
 is an Eisenstein-Jacobi integer since $A^2+3B^2=(A+B)^2-(A+B)(2B)+(2B)^2$.
 Conversely if $n$ is even then
 $m^2-mn+n^2=(m-n/2)^2+3(n/2)^2$ and since $m^2-mn+n^2=(n-m)^2-(n-m)n+n^2$ we see that every
Eisenstein-Jacobi integer is of the form $A^2+3B^2$ for some $A,B\in
\Z$. Using Proposition~\ref{ejtypet} we can write
$$A^2+3B^2=2^{2\alpha}(\underset{t\in T}{\prod}
p_t)^2\underset{j\in J} {\prod} p_j,\ \ \alpha\in \N, $$

and

$$2q=2^{2\beta}(\underset{t\in I'}{\prod} p_t)^2\underset{j\in J'}
{\prod }p_j,\ \ 1\le \beta \le \alpha, \ T' \subset T,\ J'\subset
J,$$

\n with $p_t$ prime of the form $6k-1$ for $t\in T$ and $p_j$ prime
of the form $6k+1$ or equal to $3$ for all $j\in J$. From Euler's
$6k+1$ theorem, for each $j\in J$ we can write
$p_j=(m_j+n_j\sqrt{3}i)(m_j-n_j\sqrt{3}i)$ and we make the choice of
$m_j$ and $n_j$ in $\Z$ in such a way $A+B\sqrt{3}i=s\underset{j\in
J}{\prod }(m_j+n_j\sqrt{3}i)$ using the prime factorization in
$\Z[\sqrt{3}i]$ of $A+B\sqrt{3}i$ and $h=2^{\alpha}\underset{t\in
T}{\prod} p_t$.

Then we take $r$ and $s$ such that $s+r\sqrt{3}i=h'\underset{j\in
J'}{\prod }(m_j-n_j\sqrt{3}i)$ with $h'=2^{\beta}\underset{t\in
T'}{\prod} p_t$. Notice that
$2q=(s+r\sqrt{3}i)((s-r\sqrt{3}i)=s^2+3r^2$ and
$(A+B\sqrt{3}i)(s+r\sqrt{3}i)=2q(A'+B'\sqrt{3}i)$. Identifying the
coefficients we get $As+3Br=2qA'$ and $Ar+Bs=2qB'$ and the
conclusion of our lemma follows from this.\eproof

To return to the proof of Theorem~\ref{generalpar} we begin with the
next proposition.

\begin{proposition} \label{integer} For some particular solution $(r,s)\in \Z^2$ of
 (\ref{rs}) all of the coefficients $m_u$, $m_v$, $m_w$, $n_u$, $n_v$, $n_w$, $m_x$,
 $m_y$, $m_z$, $n_x$,$n_y$, $n_z$ in (\ref{paramtwo}) are integers. \end{proposition}

\proof. \ One can check that

\begin{equation}\label{equationssbparf}
\begin{cases}
m_x^2+m_y^2+m_z^2=n_x^2+n_y^2+n_z^2=m_u^2+m_v^2+m_w^2=n_u^2+n_v^2+n_w^2=2d^2, \\
m_xn_x+m_yn_y+m_zn_z=m_un_u+m_vn_v+m_wn_w=d^2,\\
am_x+bm_y+cm_z=am_u+bm_v+cm_w=an_x+bn_y+cn_z=an_u+bn_v+cn_w=0
\end{cases}
\end{equation}

 These identities insures that the points $P(u,v,z)$ and $Q(x,y,z)$
are in the plane of normal vector $(a,b,c)$ and containing the
origin, the $\triangle OPQ$ is equilateral for every values of $m$,
$n$ and its side lengths are $l=d\sqrt{2(m^2-mn+n^2)}$. These
calculations are tedious and so we are not going to include them
here. The only ingredients that are used in establishing all these
identities are the two relations between $a$, $b$, $c$, $d$, $r$ and
$s$.

From (\ref{rs}) we see that $r$ and $s$ have to be of the same
parity. Then, it is clear that $m_z$, $n_z$, $m_w$, and $n_w$ are
all integers. Because the equalities in (\ref{equationssbparf}) are
satisfied it suffices to show that $m_x$, $n_x$, $m_u$, and $n_u$
are integers for some choice of $(r,s)$ solution of (\ref{rs}). To
show that $n_x$ is an integer we need to show that $q$ divides
$N=rac+dbs$.

Let us observe that $c^2\equiv 3d^2$ and  $a^2\equiv -b^2$ (mod
$q$). Multiplying together these two congruencies we obtain
$(ac)^2+3(db)^2\equiv 0 $ (mod $q$).  Using Lemma~\ref{existence} we
see that $N$ is divisible by $2q$ for some choice of $r$ and $s$ as
in (\ref{rs}). So, $n_x\in \Z$.

Next we want to show that $m_x$ is an integer. First let us observe
that if $M=3dbr-acs$ we can apply the second part of
Lemma~\ref{existence} to conclude that $2q$ divides $M$ also. Hence
$2q$ divides $M+N=db(3r+s)+ac(r-s)$ and so, $m_x$ is an integer.
Because $m_x+n_u=n_x$ and $m_u=n_x$ it follows that $n_u$ and $m_u$
are also integers. \eproof

{\bf Remark:} Let us observe that the parametric formulae
(\ref{paramone}) and (\ref{paramtwo}) exist under no extra
assumption on $a$, $b$ and $c$. The questions is if every triangle
in $\cal T_{a,b,c}$ is given by these formulae.

 So, with these preparations we can return to prove this under
 Theorem~\ref{generalpar}'s  assumption.

\proof. \ We start with a triangle in $\cal T_{a,b,c}$ say
$\triangle OPQ$ with the notation as before. We know that
$P(u_0,v_0,w_0)$ and $Q(x_0,y_0,z_0)$ belong to the plane of
equation $a\alpha+b\beta+c\gamma=0$ and by
Theorem~\ref{construction} we see that the coordinates of $P$ and
$Q$ should satisfy (\ref{eqfortp}). Hence, using the same notation,
$cv_0-bw_0$, $aw_0-cu_0$ and $bu_0-aw_0$ are divisible by $d$. A
relatively simple calculation shows that $m_vn_w-m_wn_v=ad$,
$m_wn_u-m_un_w=bd$ and $m_un_v-m_vn_u=cd$. We would like to solve
the following system in $m$ and $n$:

\begin{equation}\label{fsystem}
\begin{cases}
u_0=m_um-n_un,\\
v_0=m_vm-n_vn,\\
w_0=m_wm-n_wn. \\
\end{cases}
\end{equation}

\n  The equalities (\ref{equationssbparf}) and the fact that
$(u_0,v_0,w_0)$ is in the plane $a\alpha+b\beta+c\gamma=0$ insures
that (\ref{fsystem}) has a unique real solution in $m$ and $n$. We
want to show that this solution is in fact an integer solution. The
value of $n$ can be solved from each pair of these equations to get
$$n=\frac{v_0m_w-w_0m_v}{ad}=\frac{w_0m_u-u_0m_w}{bd}=\frac{u_0m_v-v_0m_u}{cd}.$$
Since $gcd(a,b,c)=1$ we can find integers $a'$, $b'$, $c'$ such that
$aa'+bb'+cc'=1$. Hence the above sequence of equalities gives
$$n=
\frac{a'(v_0m_w-w_0m_v)+b'(w_0m_u-u_0m_w)+c'(u_0m_v-v_0m_u)}{d}$$
So, it suffices to show that $d$ divides $v_0m_w-w_0m_v$,
$w_0m_u-u_0m_w$ and $u_0m_v-v_0m_u$.

Next, let us calculate for instance $v_0m_w-w_0m_v$ in more detail:
 $$\begin{array}{l}
 \ds v_0m_w-w_0m_v=v_0r-\frac{das-rbc}{q}w_0=\frac{v_0qr-(das-rbc)w_0}{q}=\\ \\
 \ds \frac{-dasw_0+v_0(3d^2-c^2)r+rbcw_0}{q}=\frac{c(bw_0-v_0c)r+3rv_0d^2-dasw_0}{q}.
\end{array}$$
From Theorem~\ref{planeprop} we see that $bw_0-v_0c=\pm
d(2x_0-u_0)$. Hence,
$$\ds v_0m_w-w_0m_v=\frac{d[\pm c(2x_0-u_0)r+3rv_0d-asw_0]}{q}.$$

This shows that $d$ divides $v_0m_w-w_0m_v$ provided that
$gcd(q,d)=1$. This is true since $gcd(d,c)=1$ implies $gcd(q,d)=1$.
Similar calculations show that $d$ divides $w_0m_u-u_0m_w$ and
$u_0m_v-v_0m_u$. Hence $n$ must be an integer. Similarly one shows
that $m$ is an integer. By replacing $P$ with $Q$ if necessary all
the equalities in (\ref{paramone}) have to hold true by
Theorem~\ref{planeprop}.\eproof

It is natural to ask if the Diophantine equation $a^2+b^2+c^2=3d^2$
has any solution which satisfies $gcd(a,b,c)=1$, $gcd(a,d)>1$,
$gcd(b,d)>1$ and $gcd(c,d)>1$. Certainly the parametrizations
(\ref{paramone}) and (\ref{paramtwo}) are the most general that one
can hope for in case the answer to the above question is negative.
Instead of  settling this question maybe one can say how rare a
solution like that can be and find in the process an asymptotic
formula for $\cal E\cal T(n)$. Let us use a similar notation here
for the number of regular tetrahedra with coordinates in the set
$\{0,1,2,...,n\}$: $\cal R\cal T(n)$. Is there any simple relation
between $\cal E\cal T(n)$ and $\cal R\cal T(n)$? What is an
asymptotic formula for $\cal R\cal T(n)$?

\end{document}